\documentclass[12pt,a4paper]{article}
\usepackage{theorem}
\usepackage{amssymb}
\usepackage{latexsym}
\usepackage{lscape}

\newtheorem{prop}{Proposition}
\newtheorem{exa}{Example}

\newtheorem{defi}{Definition}
\newtheorem{conj}{Conjecture}

\newcommand{\proof}{ \noindent{\it Proof:}\quad }
\newcommand{\qed}{\hfill$\Box$}
\makeatother 
\begin{document}
\title{\sc The Finite Field Kakeya Problem}
\author{Aart~Blokhuis and Francesco~Mazzocca}
\date{}
\maketitle \abstract{A Besicovitch set in $AG(n,q)$ is a set of
points containing a line in every direction. The Kakeya problem is
to determine the minimal size of such a set. We solve the Kakeya
problem in the plane, and substantially improve the known bounds
for $n>4$}.
\section{Introduction}
  We denote by $\pi_q$ the projective plane $PG(2,q)$ over
the Galois field $GF(q)$ with $q$ elements, $q>2$ a prime
power.
  \vskip.4cm
Let $\ell$ be a line in $\pi_q $ and, for every point
$P$ on $\ell,$ let $\ell_P$ be a line on $P$ other than $\ell .$
The set
\begin{equation}
K=(\bigcup_{P\in \ell}\ell_P)\setminus \ell
\label{def1}\end{equation} is called a \it Kakeya set\rm , or a
\it minimal Besicovitch set. \rm The \it finite plane Kakeya
problem \rm asks for the smallest size $k(q)$ of a Kakeya set; it
is the two-dimensional version of the \it finite field Kakeya
problem \rm posed by \it T.Wolff \rm in his influential paper
\cite{W} of 1996.
  \vskip.4cm In the following, unless explicitly mentioned otherwise, we will
use the same notation of (\ref{def1}) for the lines defining a
Kakeya set $K.$
  \vskip.4cm Let $\Omega$ be a set of $q+2$ points in $\pi_q.$ A
point $P\in\Omega$ is said to be an \it internal nucleus \rm of
$\Omega$ if every line through $P$ meets $\Omega$ in exactly one other point.
Internal nuclei of $(q+2)-$sets were first considered by \it
A.Bichara \rm and \it G.Korchm\'aros \rm in \cite{BK}; here they
proved the following result.
\begin{prop}\label{BK}{\rm\bf (1982)}
Let $q$ be an odd prime-power. Every set of $q+2$ points in $\pi_q$ has at most
two internal nuclei.
\end{prop}
  The $q+2$ lines  defining a Kakeya set in $\pi_q$ can be viewed as
a set of $q+2$ points with an internal nucleus in the dual plane
$\pi_q^*.$ More precisely, if $K$ is a Kakeya set in $\pi_q\, ,$
the lines $\ell$ and ${\ell}_P,$ $P\in\ell ,$ give rise in
$\pi_q^*$ to a set $\Omega (K)$ of $q+2$ points with $\ell$ as an
internal nucleus. Vice versa, every set of $q+2$ points with an
internal nucleus in $\pi_q$ defines in an obvious way a Kakeya set
in $\pi_q^*.$ Thanks to this duality, the finite plane Kakeya
problem is equivalent to ask for \it the smallest number $k^*(q)$
of lines in $\pi_q$ meeting a set of $q+2$ points with an internal
nucleus\rm ; to be precise, we have
\[
k^*(q)=1+q+k(q) \,. \label{kk} \]

\section{Old and New Results in the Plane}
  Let us start  by recalling that the first author \rm and \it
A.A.Bruen \rm studied in \cite{BB} the smal\-lest number of lines
intersecting a set of $q+2$ points in $\pi_q;$ here no
assum\-ption on the existence of internal nuclei is made.
Nevertheless the dual of the theorem $1.3$ of \cite{BB} contains
the following result as a special case.
\begin{prop}\label{BB1} {\rm\bf (1989)}
If $q\ge 7$ is odd, then
\[
|K|\ge \frac{q(q+1)}{2}+\frac{q+2}{3} \, ,
\] for every Kakeya set $K.$
\end{prop}
  \begin{exa}\label{exx1}\rm
Assume $q$ is even and consider in $\pi_q$ a dual hyperoval $\cal
H,$ i.e. a $(q+2)-$set of lines, no three of which are concurrent. Fix
a line $\ell\in\cal H$ and, for every point $P\in \ell$, let
$\ell_P$ the line of $\cal H$ on $P$ other than $\ell .$ Then the
Kakeya set
\[
K({\cal H} ,\ell )=(\bigcup_{P\in \ell }\ell_P)\setminus \ell
\] is said to be associated to $\cal H$ and
$\ell$ and it is of size
\[
|K({\cal H} ,\ell )| = \frac{q(q+1)}{2} \, . \] \qed
\end{exa}
  \begin{exa}\label{exx2}\rm
Assume $q$ is odd and consider in $\pi_q$ a dual oval $\cal O$,
i.e. a $(q+1)-$set of  lines, no three concurrent. Let
$\ell $ be a fixed line in $\cal O.$ Every point $P$ on $\ell ,$
but one, belongs to a second line $\ell_P\in \cal O$ other than
$\ell .$ If $A$ is this remaining point on $\ell$, let $\ell_A$ be a(ny)
line through it different from $\ell .$ Then the Kakeya set
\[
K({\cal O} ,\ell ,\ell_A)=(\bigcup_{P\in \ell }\ell_P )\setminus
\ell \] is said to be associated to $\cal H,$ $\ell$ and $\ell_A;$
moreover it is of size
\[
|K({\cal O} ,\ell ,\ell_A)| = \frac{q(q+1)}{2}+\frac{q-1}{2} \,
.\] \qed
\end{exa}
  For any point $A$ of a Kakeya set $K,$ we denote by $m_A$ the
number of lines $\ell_P\, ,$ $P\in\ell ,$ on $A$ and we set
\begin{equation}
\sigma (K) = \sum_{A\in K}\frac{(m_A-1)(m_A-2)}{2}\, .
\label{sum1}\end{equation} In \cite{F}, \it X.W.C.Faber \rm
described special cases of \it Examples \rm \ref{exx1} and \ref{exx2}
and, by a counting argument, proved the following result.
  \begin{prop}\label{incf}{\rm\bf (Incidence formula, 2006)} The size of a Kakeya set
$K$ is given by
\begin{equation}
|K|= \frac{q(q+1)}{2} + \sigma (K)\, . \label{sum2}\end{equation}
\end{prop}
Since $\sigma (K)\ge 0\, ,$ for every Kakeya set $K,$ a first
consequence of (\ref{sum2}) is that
\begin{equation}
|K|\ge \frac{q(q+1)}{2} \, . \label{sum3}\end{equation}  Let us
note that \it T.Wolff \rm in \cite{W} proved that $|K|\ge q^2/2;$
in fact his method gives inequality (\ref{sum3}). Equality
in (\ref{sum3}) is actually attained in \it Example \rm
\ref{exx1} and it is easy to see that this happens only in this
case. So, when $q$ is even, our problem is quite simple: \it every
Kakeya set $K$ in $\pi_q,$ $q$ even, satisfies inequality {\rm
(\ref{sum3})} and equality holds iff $K$ is associated to a dual
hyperoval and one of its lines. \rm
  \vskip.4cm When $q$ is odd the plane $\pi_q$ contains no
hyperovals and $\sigma (K)>0\, ,$ for every Kakeya set $K.$
In this case the Kakeya set closest to that of \it Example
\rm \ref{exx1} is the set $K({\cal O} ,\ell ,\ell_A)$ described in
\it Example \rm \ref{exx2}. This is the reason for the
following conjecture recently raised and studied by \it
X.W.C.Faber \rm in \cite{F}.
  \begin{conj} {\rm\bf (2006)}
If $q$ is odd, then
\[
|K|\ge \frac{q(q+1)}{2}+\frac{q-1}{2} \, , \label{fab1}\] for
every Kakeya set $K.$
\end{conj}
  \noindent We remark that the \it Blokhuis-Bruen \rm inequality in
\it Proposition \rm \ref{BB1} is not so far from that of the
conjecture. Moreover in \cite{F}, \it X.W.C.Faber \rm obtained the
following two results; the second one is a slight improvement of
\it Proposition \rm \ref{BB1}.
  \begin{prop}\label{Fa1}{\rm\bf (Triple point lemma, 2006)}
Let $K$ be a Kakeya set in $\pi_q\, ,$ $q$ odd. Then, for every point $P\in \ell$, except possibly one,
there exists a
point $A\in \ell_P$ with $m_A\ge 3.$
\end{prop}
  \begin{prop}\label{Fa2}{\rm\bf (2006)}
If $q$ is odd, then
\begin{equation}
|K|\ge \frac{q(q+1)}{2}+\frac{q}{3} \, ,
\label{fab1}\end{equation} for every Kakeya set $K.$
\end{prop}
The \it triple point lemma \rm is the the main tool in the proof
of \it Proposition \rm \ref{Fa2} and it is worth to remark that it
is just the dual of \it Proposition \rm \ref{BK}. Actually it was
proved by the same argument of \it Bichara \rm and \it
Korchm\'aros \rm : the celebrated \it Segre's lemma of tangents,
\rm that was the key ingredient in his famous cha\-racterization
of the $q+1$ rational points of an irreducible conic in $\pi_q$
with $q$ odd (\cite{Se}).
  \vskip.4cm Let $\Omega$ be a $(q+2)-$set in $\pi_q$ with an
internal nucleus  and let $\ell_{\infty}$ a line through this nucleus. Then,
in the affine plane $AG(2,q)=\pi_q\setminus\ell_{\infty},$ the
point set $\Omega\setminus\ell_{\infty}$ can be arranged as the
graph $\{ (a,f(a)) \,\, :\,\, a\in GF(q)\}$ of a function $f,$ $f$
being either a permutation or a semipermutation (i.e. a function
whose range has size $q-1$) of $GF(q).$ This graph has been
recently introduced and studied by \it J.Cooper \rm in \cite{C}
and the following improvement to the \it Faber's \rm inequality
(\ref{fab1}) has been obtained.
\begin{prop}\label{Co1}{\rm\bf (2006)}
If $q$ is odd, then
\begin{equation}
|K|\ge \frac{q(q+1)}{2}+\frac{5q}{14}-\frac{1}{14} \, ,
\end{equation} for every Kakeya set $K.$
\end{prop}

Finally, we can settle Faber's conjecture, also characterizing
the unique example realizing it. Actually we have the following
sharp result.
  \begin{prop}\label{bm1}
If $q$ is odd, then
\[
|K|\ge \frac{q(q+1)}{2}+\frac{q-1}{2} \, , \] for every
Kakeya set $K.$ Equality holds if and only if $K$ is of
type $K({\cal O} ,\ell ,\ell_A),$ as in \it Example \rm
\ref{exx2}.
\end{prop}

The essential ingredients in the proof are the \it Segre's \rm
lemma of tangents  and the \it Jamison-Brouwer-Schrijver \rm bound
on the size of blocking sets in desarguesian affine planes
(\cite{BS},\cite{J}).

\section{Solution of Kakeya's problem in the plane}

We will give the proof of Proposition \ref{bm1}. It is more convenient
however to phrase it in its dual form.

\begin{prop}\label{**}
Let $\Omega$ be a set of $q+2$ points in PG$(2,q)$, with an internal
nucleus. Then the number of lines intersecting $\Omega$ is at least
\[
k^*(q)={(q+1)(q+2)\over 2}+{q-1\over 2}.
\]
Equality implies that $\Omega$ consists of the points of an irreducible
conic together with an external point.
\end{prop}

\proof Let
$a_i$ be the number of lines in AG(2,$q$) intersecting $\Omega$ in $i$ points. Then:
\[
\left\{
\begin{array}{rcl}
\sum a_i & = & q^2+q+1\\
\sum i\,a_i & = & (q+2)(q+1)\\
\sum {i\choose 2}\,a_i &=& (q+2)(q+1)/2
\end{array}
\right.
\]
The first equation counts the total number of lines in the affine plane.
In the second we count incident point-line pairs $(P,\ell)$, where $P$ is a point of
$\Omega$. Finally in the third we count ordered triples $(P,Q,\ell)$, where
$P$ and $Q$ are different points from $\Omega$ (and $\ell$ the unique line joining
them). It follows that
\[
a_0+a_3+3a_4+\dots+{q\choose 2}a_{q+1}=(q^2-q)/2.
\]
Also, for later use we note that:
\[
a_1=3a_3+8a_4+\dots=\sum_{n>2} (n^2-2n)a_n.
\]
We aim for the situation where $\Omega$ is a conic together with an external point.
In that case $a_1=(q-1)+(q-1)/2$, $a_2=(q^2+5)/2$, $a_3=(q-1)/2$
and $a_0=(q-1)^2/2$ (and the number of intersecting lines is $(q^2+4q+1)/2$).\\
Let the number of intersecting
lines be $(q+2)(q+1)/2+f$ for some $f$, so that $a_0=(q^2-q)/2-f$. This gives us for $f$
the equation
\[
a_3+3a_4+\dots+{q\choose 2}a_{q+1}=f,
\]
and we would like to show that $f\ge (q-1)/2$.\\

We know from {\it Bichara-Korchm\'aros} result (Prop.\ref{BK}),
that there are at most 2 internal nuclei (in the example exactly
2) and by assumption there is at least one. Every other point is
therefore on at least one tangent, and hence also on at least one
$(\ge3)$-secant. In particular $f\ge q/3$, with equality if every
other point is on exactly one tangent and
one three-secant (this does happen if $q=3$).\\

Every point, with the exception of the internal nucleus (nuclei), is on an odd intersector.
So the odd intersectors form a blocking set of the dual affine plane if there is just
one nucleus (this should maybe be called a dual blocking set, but we will use this
term with a different meaning later). In this case:
\[
a_1+a_3+a_5+\dots \ge 2q-1,
\]
and therefore
\[
4a_3+8a_4+15a_5+\dots \ge 2q-1,
\]
and hence $f\ge (2q-1)/4$, more than we want.\\
From now on we assume that there are two internal nuclei, $N_1$ and $N_2$. Adding
a random line on one of the internal nuclei, but not containing the other one, we again
get a blocking set of the dual affine plane,
and we obtain
\[
4a_3+8a_4+15a_5+\dots \ge 2q-2,
\]
and hence $f\ge (2q-2)/4$ with equality if $a_k=0$ for $k>3$. So we have proved our lower bound,
and we proceed to characterize the case of equality.\\
If $f=(q-1)/2$ then we have $(q-1)/2$ three-secants, and $3(q-1)/2$ tangents. Now if a point
$Q$, is on exactly one tangent, and this happens often, then also on a unique three-secant,
and  we will show, that their intersection points with $\ell$ are related: if one is $(1:\lambda)$ the
other is $(1:-\lambda)$), where coordinates are chosen such that $N_1=(1:0)$ and $N_2=(0:1)$.\\

Consider a three-secant containing two points on a unique tangent.
Then these two tangents intersect in a point on the line joining
the two internal nuclei ($\ell$). This is true in the example and
follows from a {\it Segre}-type computation: if the three secant
intersects the line $\ell$ in $(1:\lambda:0)$ then the unique
tangents go through $(1:-\lambda:0)$), where the coordinates are
set up in such a way that the two
internal nuclei are $(1:0:0)$ and $(0:1:0)$.\\
We will use {\it Segre}-type computations a lot in the sequel. The
general setup is the following. Consider three points
$E_1=(1:0:0)$, $E_2=(0:1:0)$, $E_3=(0:0:1)$. Let $X$ be any set of
points such that no point of $X$ is on one of the coordinate lines
$E_iE_j$. For $x=(x_1:x_2:x_3)$ write down the triple
$x'=(x'_1,x'_2,x'_3):=(x_2/x_1, x_3/x_2, x_1/x_3)$. It is clear
from the definition that $\prod_{x\in X} x'_1x'_2x'_3=1$. On the
other hand, it is sometimes possible, because of geometric
properties of $X$ to say something about $p_i=\prod_{x\in X}x'_i$.
Applying this together with $p_1p_2p_3=1$ is called {\it Segre}'s
lemma of tangents or a {\it Segre} computation.  In our case the
argument runs as follows. Let $U$ be a point on a unique
three-secant, further choose coordinates such that $U=(0:0:1)$,
and some random fourth point equals $(1:1:1)$. Recall that
$N_1=(1:0:0)$ and $N_2=(0:1:0)$. Let the three-secant through $U$
intersect $\ell$ in $(1:\lambda:0)$ and let the unique tangent
intersect $\ell$ in $(1:\mu:0)$. The remaining $q-1$ points of
$\Omega$ (other than $N_1$, $N_2$ and $U$) have (homogeneous)
coordinates $(a_i:b_i:c_i)$ with $a_ib_ic_i\ne0$. We associate to
such a point the triple $(b_i/a_i,c_i/b_i,a_i/c_i)$. Taking the
product of all the entries in all triples we clearly get 1,
because that is the contribution of each triple. On the other hand
we have $\prod_i c_i/b_i=-1$, because on each line through $N_1$
we have a unique point of $\Omega$ so we just have the product of
all non-zero field elements. In the same way $\prod a_i/c_i=-1$ by
considering lines through $N_2$. To compute $\prod b_i/a_i$ we
consider the lines through $U=(0:0:1)$. The three secant gives the
value $b_i/a_i=\lambda$ twice, but the value $b_i/a_i=\mu$ is
absent. All other nonzero field elements occur exactly once in the
product, so for this product we end up with $-\lambda/\mu$, so
$(-1)(-1)(-\lambda/\mu)=1$ and we
conclude that $\mu=-\lambda$.\\

We will show that, unless $q=3$, the three points of $\Omega$ on a three-secant cannot all be points
with a unique tangent, by applying again a {\it Segre} computation.\\
Apart from the 2 internal nuclei our set has $q$ points, and all
of them are on at least one tangent. The total number of tangents
is $$3(q-1)/2=q+(q-3)/2$$ hence at least $(q+3)/2$ points are on
exactly one tangent (and one three-secant). So we certainly find a
three-secant
with (at least) two unique-tangent points on it.\\
Let $N_1=(1:0:0)$ and $N_2=(0:1:0)$ (as before) be the internal nuclei.\\
Let $U_1=(0:0:1)$ and $U_2=(1:1:1)$ be two one-tangent points on a
common three-secant, and let $V=(a:b:1)$ be
a one-tangent point not on the line $U_1U_2$, so $a$ and $b$ are nonzero, and $a\ne b$.\\
Note that in our example we have that $N_1, N_2, U_1$ and $U_2$
are on a conic, and the tangents at $U_{1,2}$ are also known. So
the conic has to be: $-2x_1x_2+x_2x_3+x_3x_1=0$. So we should
expect that $-2ab+a+b=0$ for $V=(a:b:1).$ \\
The three-secant $U_1U_2$ meets $N_1N_2$ in $(1:1:0)=N_1+N_2$, so
the tangents at $U_1$ and $U_2$ meet in $N_1-N_2=(1:-1:0)$. Let
the tangent at $V$ pass through $(1:\lambda:0)$, then the
three-secant on  $V$ passes through
the point $(1:-\lambda:0)$.\\
First we consider the triangle $U_1N_1V$. The tangent at $V_1$
intersects $U_1N_1$ in $U_1+((\lambda a-b)/\lambda)N_1$, the
three-line in $U_1+((\lambda a+b)/\lambda)N_1$. The tangent
through $U_1$ intersects $N_1V$ in $N_1+(-1/(a+b))V$, the
three-line in $N_1+(1/(b-a))V$. On $VU_1$ there are no special
'missing' or 'extra' points. {\it Segre} gives:
\[
(a+b)(\lambda a+b)=(b-a)(\lambda a-b).
\]
And we get the important fact $\lambda=-b^2/a^2$. \\
Next we consider the triangle $N_1U_2U_1$. Let the third point of $\Omega$ on
$U_1U_2$ be $U_2+\mu U_1$. On $N_1U_2$ we 'miss' the point $(-1:1:1)=N_1+(-1/2)U_2$. On
$U_2U_1$ we 'miss' the point $U_2+\mu U_1$, and finally on $U_1N_1$ the point $(2:0:1)=U_1+2N_1$.
Here we used that since the three-line on $U_1$ goes through $(1:1:0)$, the tangent passes through
$(1:-1:0)$. It follows from the {\it Segre} product that $\mu=1$. \\
We now turn to
the triangle $U_1U_2V_1$. On $U_1U_2$ we find the 'extra' point, the intersection with the
three line through $V$:
\[
U_1+{(b+a\lambda)/(1+\lambda)\over 1-(b+a\lambda)/(1+\lambda)}U_2.
\]
and 'missing' points $U_1+U_2$ (the third point of $\Omega$ on $U_1U_2$) and the intersection of
the tangent through $V$ with $U_1U_2$:
\[
U_1+{(b-a\lambda)/(1-\lambda)\over 1-(b-a\lambda)/(1-\lambda)}U_2.
\]
This is of course just the expression for the three-secant with
$-\lambda$ instead of $\lambda$. On $U_2V$ and $VU_1$ we find
'missing' coordinates $-2/(a+b)$ and $-1+(a+b)/2$. The {\it Segre}
computation gives us
\[
(a+b)(b+a\lambda) (1-b+(a-1)\lambda)=(a+b-2)(b-a\lambda)(1-b-(a-1)\lambda).
\]
This we may rewrite as
\[
a(a-1)\lambda^2+(a+b-1)(a-b)\lambda +b(1-b)=0.
\]
Now substitute $\lambda=-b^2/a^2$, multiply by $a^3$ and divide by $b$. We get:
\[
(a-b)(a+b)(2ab-a-b)=0.
\]
We already remarked that $a\ne b$, but also $a\ne -b$ because otherwise $V$ would be on the
tangent through $U_1$. Hence $2ab-a-b=0$ and $V$ is a point on the conic we are aiming for.
A direct computation shows that also the tangent is 'right' and that the three-secant through
$V$ passes through the 'special point' $(1:1:2)=U_1+U_2$.\\
Some counting to end the story. Let there be $k$ points on a unique tangent. This means that
our special point $U_1+U_2$ is on at least $k/2$ three-secants, and hence on at least $k/2$
tangents. What is left in $\Omega$ (apart from the internal nuclei, the special point and the
unique tangent points) is a set of $q-1-k$ points on at least 2 tangents, and a set
of at most $3(q-1)/2-k-k/2$ tangents. So
\[
3(q-1)/2-k-k/2\ge 2(q-1-k).
\]
This means $k\ge q-1$, so all other points are on the conic, and we finished the proof.

\section{Applications to Dual Blocking Sets}
A blocking set $B$ in $\pi_q=PG(2,q)$ is a point set meeting every
line and containing none.
  \begin{defi}
A dual blocking set $S$ in $\pi_q$ is a point set meeting every
blocking set and containing no lines.
\end{defi}
  \begin{exa}\label{exx3}\rm
A Kakeya set $K=(\bigcup_{P\in \ell}\ell_P)\setminus \ell$  in
$\pi_q$ contains no lines. Moreover, for every blocking set $B$ of
$\pi_q,$ a point $P$ exists on $\ell\setminus B$ and so $K$ meets
$B$ in a point of $\ell_P\setminus \ell .$ It follows that $K$ \it
is a dual blocking set. \rm \qed
\end{exa}
  \begin{exa}\label{exx4}\rm
The complement $S=\pi_q\setminus (\ell \cup m)$ of the union of
two distinct lines $\ell$ and $m$ in $\pi_q$ contains no lines.
Moreover, no blocking set is contained in the union of two lines
and so $S$ meets every blocking set. It follows that $S$ \it is a
dual blocking set. \rm \qed
\end{exa}

  \vskip.4cm Dual blocking sets were introduced by \it
P.Cameron, \it F.Mazzocca \rm and \it R.Me\-shu\-lam \rm
in \cite{CMM}; the first of the two main results of this
paper is the following.
  \begin{prop}\label{cmm1} {\rm\bf (1988)}
Let $S$ be a dual blocking set in $\pi_q.$ Then
$$|S|\ge\frac{q(q+1)}{2}.$$ Equality holds if and only if
either
\begin{itemize}
  \item [\rm (i)] $S$ is the Kakeya set associated
  to a dual hyperoval and one of its lines; or
  \item [\rm (ii)] $q=3$ and $S$ is the complement
  of the union of two distinct lines.
\end{itemize}
\end{prop}
The argument in the proof of this proposition implicitly shows
that every minimal (with respect to inclusion) dual blocking set
in $\pi_q$ is of one of types described in examples (\ref{exx3})
and (\ref{exx4}). For the sake of completeness we give an explicit
proof of this result.
\begin{prop} \label{cmm2}
Let $S$ be a minimal dual blocking set in $\pi_q.$ Then one of the
two following possibilities occur:
\begin{itemize}
  \item [\rm (i)] $S=
  (\bigcup_{P\in \ell}\ell_P)\setminus \ell$ is a Kakeya
  set;
  \item [\rm (ii)] $S=\pi_q\setminus (\ell \cup m)$ is the complement of the union of
two distinct lines $\ell$ and $m.$
\end{itemize}
\end{prop}
\proof First of all we observe that there is a line
$\ell$ disjoint from $S$, for if not, then, since $S$ does not contain a
line, $S$ and its complement are blocking sets; a
contradiction as $S$ must meet every blocking set. Now we
distinguish the following two cases.\\
  \it Case 1. \rm Assume that $S$ is disjoint from exactly one
line $\ell $, and let $P$ be a point of this line. If, for
every line $m\not= \ell$ through $P$, there is a point
$Q\not= P$ on $m$ but not in $S$, then
$$ B= \left( \ell\setminus\{ P\} \right) \cup
\left(\bigcup_{P\in m\not=\ell} m  \right) $$ is a
blocking set disjoint from $S$; a contradiction. Hence,
for every point $P\in\ell ,$ there exists a line $\ell_P$
through $P$ with $\ell_P\setminus\{ P\} \subseteq S.$
Then $S$ contains the Kakeya set $K=(\bigcup_{P\in
\ell}\ell_P)\setminus \ell ,$ which is a dual blocking
set. From the minimality of $S$ it follows that $S=K$.\\
  \it Case 2. \rm Assume that there are two lines $\ell$ and  $m $
disjoint from $S$. For any point $P\not\in\ell\cup m$ ,let
$n$ be a line on $P$ meeting $\ell\setminus m$ and
$m\setminus\ell$ in the points $L$ and $M,$ respectively.
Then $(\ell\cup m\cup\{ P\} )\setminus\{ L,M\} $ is a
blocking set contained in $\ell\cup m\cup\{ P\} $. It follows that
$P$ must belong to $S$ and $S$ is the complement of
$\ell\cup m.$ \qed
  \vskip.4cm By \it Propositions \rm \ref{cmm1} and \ref{cmm2}
 we can conclude that all the bounds previously shown for
the size of a Kakeya set give, in the case that $q$ is odd,
corresponding new
bounds for the size of a minimal dual blocking set,
improving the result of \it Proposition \rm \ref{cmm1}.
In fact, as a corollary of \it Proposition \rm
\ref{bm1}, we have the following sharp result.
  \begin{prop}\label{bm3} Let $S$ be a dual blocking
set in $\pi_q,$ $q$ odd. Then
\[
|S|\ge \frac{q(q+1)}{2}+\frac{q-1}{2} \,  \]  and
equality holds if and only if $S$ is a Kakeya set of type
described in \it Example \rm \ref{exx2}.
\end{prop}

\section{Old and new results in higher dimensions}

In contrast to the plane case we only have bounds and conjectures
for higher dimensions. In \cite{W} it is shown that the number of
points in a Kakeya set in AG$(n,q)$ is at least $c\cdot
q^{(n+2)/2}$, which is good for $n=2$ but probably not for any
larger $n$. The case $n=3$ is the first open problem, but for
$n=4$ \it T.Tao \rm has shown (\cite{T}) that the exponent 3 can
be improved to $3+{1\over 16}$. In what follows we will show that
for general $n$ we get the lower bound $c\cdot q^{n-1}$, where
$c=1/(n-1)!$, so this improves the previous bounds when $n$ is at
least 5 and comes close to the conjectured $c_n q^n$.
Unfortunately our ideas are for several reasons very
unlikely to lead to improvements in the case of the 'real' Kakeya problem.\\
Very recently however, \it Zeev Dvir \rm \cite{D} has proved the
finite field Kakeya problem, by showing that the number of points
of a Kakeya set in AG$(n,q)$ is at least $q+n-1\choose n$.

Since our result and proof are similar in nature but still
slightly different, we will include it for historical reasons, and
with the hope that an improved argument will give a bound
equivalent or even slightly better than that of \it Dvir. \rm

To improve the bound in higher dimensions we use a bound
on the dimension of a certain geometric codes.\\

Consider the line-point incidence
matrix of $PG(n,q)$. Number the points (so the columns): first the points in
the hyperplane at
infinity, then the points not in the Kakeya set, and finally the points in
the Kakeya set.
As usual we denote the number of points (and hyperplanes)
in $PG(n,q)$ by $\theta_n=\theta_n(q)=(q^{n+1}-1)/(q-1)$.\\
Let the first $\theta_{n-1}$ rows be labeled by the lines
defining the Kakeya set, in the right order.
The top consisting of the first $\theta_{n-1}$ rows of the incidence matrix now looks like
this:
\[
{\bf T}=
\left(
{\bf I}\,;\,{\bf O}\,;\,{\bf K}
\right).
\]
Here ${\bf I}$ is the identity matrix, and ${\bf K}$ is the $\theta_{n-1}$ by $|K|$
line-point incidence matrix of Kakeya-lines versus Kakeya-points
Let $d=d_{n-1}$ be the dimension of $C_{n-1}$, the $GF(p)$-code (where $q=p^t$) spanned by the
lines of $PG(n-1,q)$ (the hyperplane at infinity). Then there is a subset $C$
of the points, of size
$\theta_{n-1}-d_{n-1}$ that does not contain the support of
a codeword (this is obvious: after normalization a generator matrix for this
code has the form $(I\,;\,A)$ and every nonzero codeword has a nonzero coordinate in one
of the
first $d_{n-1}$ positions, so no codeword has its support contained in the 'tail' of
length $\theta_{n-1}-d_{n-1}$). It follows that the set of Kakeya points has at
least this size: Consider the $\theta_{n-1}-d_{n-1}$ rows of ${\bf T}$ corresponding to the
Kakeya lines having a direction in $C$. Suppose the corresponding rows of ${\bf K}$ are
dependent (over $GF(p)$).
Then this dependency would produce a codeword in the line-point code of
PG$(n,q)$ with support contained in the set $C$ in the hyperplane at infinity. But such
a word is already in the point-line code of this hyperplane. To see this, let $C_n$ stand
for the line code of $PG(n,q)$, and $C_{n-1}$ for the line code of the hyperplane $H$.
Clearly $C_n^\perp|_H\subseteq C_{n-1}^\perp$. We show that in fact equality holds, for
let $u$ be a word in $C_{n-1}^\perp$, and now take a point $P\not\in H$ and form the
cone with top $P$ over $u$, but remove $P$. This defines in an obvious way a word $\tilde u$
in $C_n^\perp$ whose restriction to $H$ is $u$.\\
So we find $|K|\ge \dim C_{n-1}^\perp$.
The dimension of $C_{n-1}$ is known,
and equal to something complicated. For us the bound
\[
|K|\ge \dim C_{n-1}^\perp \ge {q+n-2\choose n-1}\ge q^{n-1}/(n-1)!
\]
suffices. In fact, if $q$ is prime we have equality, if not we have a little improvement,
but not an essential one.

\vspace{10pt} \noindent {\bf ACKNOWLEDGMENTS } This research was
done when the first author was visiting the \it Seconda
Universit\'a degli Studi di Napoli \rm at Caserta. The authors
wish to thank for their supports the research group GNSAGA of
Italian \it Istituto Nazionale di Alta Ma\-te\-ma\-ti\-ca \rm and
the \it Mathema\-tics Department \rm of the \it Seconda
Universit$\grave{a}$ degli Studi di Napoli. \rm

\vspace{0.3cm}

\noindent{\em Authors address:}

{\sc A.Blokhuis} \par

Eindhoven University of Technology \par

Department of Mathematics and Computing Science \par

P.O. Box 513 \par

5600 MB Eindhoven, The Nederlands \par

e-mail: \verb"a.blokhuis@tue.nl"

\vspace{0.2cm}

{\sc F.Mazzocca} \par

Seconda Universit\`{a} degli Studi di Napoli\par

Dipartimento di Matematica \par

via Vivaldi 43, \par I-81100 Caserta, Italy \par

e-mail: \verb"francesco.mazzocca@unina2.it"

\end{document}